\documentclass[11pt]{article}
\usepackage{amssymb}

\usepackage{amsmath}

\input{tcilatex}

\begin{document}

\title{On the inverse Kostka matrix}
\author{Haibao Duan \\
Institute of Mathematics, Chinese Academy of Sciences, \\
Beijing 100080, dhb@math.ac.cn}
\maketitle

\begin{abstract}
In the ring of symmetric functions the inverse Kostka matrix appears as the
transition matrix from the bases given by monomial symmetric functions to
the Schur bases.

We present both a combinatorial characterization and a recurrent formula for
the entries of the inverse Kostka matrix which are different from the
results obtained by Egecioglu and Remmel [ER] in 1990. An application to the
topology of the classifying space BU(n) is obtained.
\end{abstract}



\section{Introduction}

By a \textsl{partition }we mean a sequence $\lambda =(\lambda _{1},\ldots
,\lambda _{n})$ of $n$ non-negative integers in non-decreasing order $0\leq
\lambda _{1}\leq \ldots \leq \lambda _{n}$. The following alternative
notions for partitions will be useful in simplifying our presentation.

1) For a sequence $0\leq \lambda _{1}\leq \ldots \leq \lambda _{k}$ of $k$
integers with $k\leq n$, $\lambda =(\lambda _{1},\ldots ,\lambda _{k})$
stands for the partition that differs from $\lambda $ only by $n-k$ of
zero's at the beginning.

2) If $\lambda _{1}\geq 1$ we may write $\lambda =\{\lambda _{1},\ldots
,\lambda _{k}\}$ instead of $\lambda =(\lambda _{1},\ldots ,\lambda _{k})$.

3) Sometimes we use $\lambda =(r_{1}^{i_{1}},\ldots ,r_{k}^{i_{k}})$, where

$\quad 1\leq r_{1}<\ldots <r_{k}$; $1\leq i_{1},\ldots ,i_{k}$

\noindent and where $\Sigma i_{s}\leq n$, to indicate that $\lambda $ is the
partition which begins with $n-\Sigma i_{s}$ copies of $0$, followed by $%
i_{1}$ copies of $r_{1}$, then $i_{2}$ copies of $r_{2}$, $\ldots $, etc.

Every partition can be uniquely written in the form of 2) (resp. of 3)). If
this is the case, the number $k$ (resp. $\Sigma i_{s}$) will be called
\textsl{the length of} $\lambda $ and will be denoted by $l(\lambda )$.

We say $\lambda $ is a partition of $m$ if $m=$ $\lambda _{1}+\ldots
+\lambda _{n}$. Denote by $P(m\mid n)$ for the set of all partitions of%
\textsl{\ }$m$. Assume throughout that $n\geq m$.

\bigskip

Let $\Lambda _{n}$ be the ring of symmetric functions in the variables $%
x_{1},\ldots ,x_{n}$. It is graded by $\Lambda _{n}=\bigoplus\limits_{m\geq
0}\Lambda _{n}^{m}$, where $\Lambda _{n}^{m}$ is the $\mathbb{Z}$-module
consisting of all homogenous symmetric polynomials of degree $m$, together
with the zero polynomial. We recall two canonical $\mathbb{Z}$-bases of $%
\Lambda _{n}^{m}$, both are parameterized by elements in $P(m\mid n)$.

\bigskip

Let $\mathbb{N}^{n}$ be the set of $n$-tuples of non-negative integers. The
set of partitions $P(m\mid n)$ will be considered as a subset of $\mathbb{N}%
^{n}$ in the obvious way. For $\alpha =(\alpha _{1},\ldots ,\alpha _{n})\in
\mathbb{N}^{n}$ write $x^{\alpha }$ for the monomial $x_{1}^{\alpha
_{1}}\ldots x_{n}^{\alpha _{n}}$.

Let $S_{n}$ be the permutation group on $\mathbb{N}^{n}$, and let $%
S_{n}(\alpha )$ be the stabilizer of $S_{n}$ at $\alpha \in \mathbb{N}^{n}$.
The set of left cosets of $S_{n}(\alpha )$ in $S_{n}$ is denoted by $%
S_{n}^{\alpha }$.

\textbf{Definition.} The \textsl{monomial symmetric functions }$m_{\lambda
}(n)$ associated to\textsl{\ }a $\lambda \in P(m\mid n)$ is the element of $%
\Lambda _{n}^{m}$ defined by $m_{\lambda }(n)=\sum\limits_{w\in
S_{n}^{\lambda }}x^{w(\lambda )}$.

For $\alpha =(\alpha _{1},\ldots ,\alpha _{n})\in \mathbb{N}^{n}$ let $%
a_{\alpha }=\det (x_{j}^{\alpha _{i}})_{n\times n}$. The \textsl{Schur
functions }$s_{\lambda }(n)$ associated to\textsl{\ }$\lambda =(\lambda
_{1},\ldots ,\lambda _{n})\in P(m\mid n)$ is the element of $\Lambda
_{n}^{m} $ defined by $s_{\lambda }(n)=a_{\lambda +\delta (n)}/a_{\delta
(n)} $, where $\delta (n)=(0,1,\ldots ,n-1)$ and where the addition $\lambda
+\delta (n)$ takes place in $\mathbb{N}^{n}$.$\square $

\bigskip

It is well known that both of the two sets of functions $\{m_{\lambda
}(n)\mid \lambda \in P(m\mid n)\}$ and $\{s_{\lambda }(n)\mid \lambda \in
P(m\mid n)\}$ constitute additive bases of $\Lambda _{n}^{m}$ (cf. [Ma]). An
immediate consequence is the existence of transition matrixes (with integer
entries) from one to the other. Namely, we have

$\qquad \ \ \ s_{\lambda }(n)=\sum\limits_{\mu \in P(m\mid n)}K_{\lambda
,\mu }m_{\mu }(n)$;

\noindent (1.1) $\ \ \qquad m_{\lambda }(n)=\sum\limits_{\mu \in P(m\mid
n)}K_{\lambda ,\mu }^{-1}s_{\mu }(n)$

\noindent for some integer matrices $K_{m}=(K_{\lambda ,\mu })$, $%
K_{m}^{-1}=(K_{\lambda ,\mu }^{-1})$ with $K_{m}K_{m}^{-1}=Id$.

\bigskip

The matrix $K_{m}$ (resp. $K_{m}^{-1}$) is known as the \textsl{Kostka matrix%
} (resp. the \textsl{inverse} \textsl{Kostka matrix}). They were introduced
by Kostka in 1882 [K1], who also computed $K_{m}$ and $K_{m}^{-1}$ up to $%
m=11$ [K2]. Despite much combinatorial interest in these matrices (cf. [Ma,
p.99-111]) we choose to emphasize a connection of a topological problem to
the question of effective computations of certain $K_{\lambda ,\mu }^{-1}$.

\bigskip

Let $BU(n)$ (resp. $BO(n)$) be the classifying space for the complex unitary
group $U(n)$ (resp. the real orthogonal group $O(n)$) of order $n$, and let $%
1+c_{1}+\ldots +c_{n}\in H^{\ast }(BU(n);\mathbb{Z})$ (resp. $1+w_{1}+\ldots
+w_{n}\in H^{\ast }(BO(n);\mathbb{Z}_{2})$) be the total Chern class of the
universal complex $n$-bundle over $BU(n)$ (resp. the total Stiefel-Whitney
class of the universal real $n$-bundle over $BO(n)$) [MS, \S 5]. Then

$H^{\ast }(BU(n);\mathbb{Z})=\mathbb{Z}[c_{1},\ldots ,c_{n}]$\quad (resp. $%
H^{\ast }(BO(n);\mathbb{Z}_{2})=\mathbb{Z}_{2}[w_{1},\ldots ,w_{n}]$).

\noindent As is classically known, the correspondence $c_{m}\rightarrow
m_{(1^{m})}$ (resp. $w_{m}\rightarrow m_{(1^{m})}$ mod $2$) establishes a
grading preserving isomorphism

$\qquad H^{\ast }(BU(n);\mathbb{Z})\cong \Lambda _{n}$ (resp. $H^{\ast
}(BO(n);\mathbb{Z}_{2})\cong \Lambda _{n}\otimes \mathbb{Z}_{2}$) [MS, \S 7].

\noindent In particular, the symmetric functions $s_{\lambda }(n)$ and $%
m_{\lambda }(n)$ (resp. module $2$) can be interpreted as closed cocycles on
$BU(n)$ (resp. on $BO(n)$) in which the $s_{\lambda }(n)$ act as the
Kronecker dual of the classical Schubert varieties $\Omega (\lambda )$ in $%
BU(n)$ (resp. in $BO(n)$) [MS, \S 6].

\noindent

Let $\mathcal{P}^{k}$ (resp. let $Sq^{k}$) be the Steenrod mod-$p$
operations for a prime $p\neq 2$ on the mod-$p$ (resp. the Steenrod mod-$2$
operations on the mod-$2$) cohomology of spaces [St]. From the Cartan
formula [St] we get

$\qquad \mathcal{P}^{k}(c_{m})=m_{(1^{m-k},p^{k})}$ mod $p$ in $H^{\ast
}(BU(n);\mathbb{Z}_{p})$

\qquad (resp. $Sq^{k}(w_{m})=m_{(1^{m-k},2^{k})\text{ }}$mod $2$ in $H^{\ast
}(BO(n);\mathbb{Z}_{2})$),

\noindent where $k\leq m$. Combining this with (1.1) yields

\noindent (1.2) $\ \ \mathcal{P}^{k}(c_{m})=\sum\limits_{\mu \in
P(m+(p-1)k\mid n)}K_{(1^{m-k},p^{k}),\mu }^{-1}s_{\mu }$ mod $p$.

\ \ \ \ (resp. $Sq^{k}(w_{m})=\sum\limits_{\mu \in P(m+k\mid
n)}K_{(1^{m-k},2^{k}),\mu }^{-1}s_{\mu }$ mod $2$).

\noindent In view of the decomposition of $BU(n)$ given by Schubert cells,
the geometric significance of (1.2) is clear: the number $%
K_{(1^{m-k},p^{k}),\mu }^{-1}$ gives information on how the cell $\Omega
(\mu )$ is attached to $\Omega (1^{m})$. We quote from Lenart [L]: very
little is known about the attaching maps of Schubert cells.

\bigskip

In an attempt to find generalization of the famous Wu-formula [W]

\noindent (1.3)$\qquad Sq^{k}(w_{m})=\sum\limits_{0\leq i\leq k}\binom{m-i-1%
}{k-i}w_{i}w_{m+k-i}$ (mod $2$),

\noindent many works were devoted to express $\mathcal{P}^{k}(c_{m})$ as a
polynomial in the Chern classes $c_{i}$ (See Shay [S] or Lenart [L] for the
history and relevant works on this). However, this task may alternatively be
implemented by combining (1.2) with the Giambelli-formula [Ma] which
expresses $s_{\mu }$ as a polynomial in the $c_{i}$.

\bigskip

Turning to the numerical aspects of the matrix $K_{m}^{-1}$, Egecioglu and
Remmel [ER] obtained in 1990 a combinatorial interpretation and a recurrent
formula for $K_{\lambda ,\mu }^{-1}$. We briefly recall their result.

\bigskip

For a partition $\lambda =(r_{1}^{i_{1}},\ldots ,r_{k}^{i_{k}})$ and $1\leq
j\leq k$, let $\lambda \lbrack j]$ denote the partition which results from $%
\lambda $ by removing a part $r_{j}$. That is $\lambda \lbrack
j]=(r_{1}^{i_{1}},\ldots
,r_{j-1}^{i_{j-1}},r_{j}^{i_{j}-1},r_{j+1}^{i_{j+1}}\ldots ,r_{k}^{i_{k}})$.

If $\mu =\{\mu _{1},\ldots ,\mu _{l}\}$ is a partition and $1\leq i\leq l$,
we write $\omega \subset _{i}\mu $ for $\omega =(\mu _{1}-1,\ldots ,\mu
_{i-1}-1,\mu _{i+1},\ldots ,\mu _{l})$.

For two partitions $\lambda =\{\lambda _{1},\ldots ,\lambda _{k}\},\mu
=\{\mu _{1},\ldots ,\mu _{l}\}\in P(m\mid n)$ denote by $T(\lambda ,\mu )$
the set of all chains of partitions

$T:\qquad \mu ^{0}=(0)\subset _{j_{1}}\mu ^{1}\subset _{j_{2}}\ldots \subset
_{j_{k}}\mu ^{k}=\mu $

\noindent so that if we have $\mu ^{i}=\{\mu _{1}^{\prime },\ldots ,\mu
_{l^{\prime }}^{\prime }\}$ and put $a_{i}=\mu _{j_{i}}^{\prime }+j_{i}-1$
for each $1\leq i\leq k$, then the sequence $(a_{1},\ldots ,a_{k})$ is a
permutation of $\lambda _{1},\ldots ,\lambda _{k}$.

For each $T\in T(\lambda ,\mu )$ we set $sign(T)=(-1)^{j_{1}+\ldots
+j_{k}-k} $\

\textbf{Theorem 1 }(Egecioglu-Remmel [ER]). \textsl{If }$\lambda
=(r_{1}^{i_{1}},\ldots ,r_{k}^{i_{k}})$\textsl{\ and }$\mu =\{\mu
_{1},\ldots ,\mu _{l}\}$\textsl{, then}

\noindent (1.4)$\qquad K_{\lambda ,\mu }^{-1}=\sum\limits_{\substack{ 1\leq
j\leq d  \\ r_{j}=\mu _{i}+i-1}}(-1)^{i-1}K_{\lambda \lbrack j],(\mu
_{1}-1,\ldots ,\mu _{i-1}-1,\mu _{i+1},\ldots ,\mu _{l})}^{-1}$.

\noindent \textsl{Consequently, }$K_{\lambda ,\mu }^{-1}=\sum\limits_{T\in
T(\lambda ,\mu )}sign(T)$\textsl{. }

\bigskip

The merit of the recurrence (1.4) is apparent: without resorting to the
symmetric functions (i.e. (1.1)), the computation of $K_{\lambda ,\mu }^{-1}$
can be boiled down to the obvious relation $K_{(0),(0)}^{-1}=1$.

It might be worthwhile to mention that Lenart announced in [L, Remark 5.4]
also an algorithm to determine $K_{\lambda ,\mu }^{-1}$.

\bigskip

In this paper we give another combinatorial description and recurrence for
the numbers $K_{\lambda ,\mu }^{-1}$. Our main results are stated in Section
2, followed by their proofs in Section 3. Our formula (2.1) is ready to
apply to yield computational results. Examples are shown in Section 4.

\section{Main results}

For a $\lambda \in P(m\mid n)$ define a map $f_{\lambda }:S_{n}^{\lambda
}\times S_{n}\rightarrow \mathbb{N}^{n}$ by letting

$\qquad \qquad f_{\lambda }(w,\sigma )=w(\lambda )+\sigma (\delta (n))$.

\noindent Let $\varepsilon :S_{n}\rightarrow \{\pm 1\}$ be the sign function
on the permutation group $S_{n}$.

\textbf{Theorem 2. }\textsl{Assume that }$m_{\lambda }(n)=\sum\limits_{\mu
\in P(m\mid n)}K_{\lambda ,\mu }^{-1}s_{\mu }(n)$\textsl{. Then}

$\qquad \qquad K_{\lambda ,\mu }^{-1}=\sum\limits_{(w,\sigma )\in f_{\lambda
}^{-1}(\mu +\delta (n))}\varepsilon (\sigma )$\textsl{.}

\textbf{Remark 1.} Theorem 2 offers another interpretation of the numbers $%
K_{\lambda ,\mu }^{-1}$. Given two partitions $\lambda ,\mu \in P(m\mid n)$,
the set $f_{\lambda }^{-1}(\mu +\delta (n))$ is of interest for it consists
of solutions $(w,\sigma )\in S_{n}^{\lambda }\times S_{n}$ to the vector
equation

$\qquad \qquad w(\lambda )+\sigma (\delta (n))=\mu +\delta (n)$.

\noindent Let $I:S_{n}\rightarrow \mathbb{N}$ be the length function on the
permutation group. From the subset $f_{\lambda }^{-1}(\mu +\delta
(n))\subset S_{n}^{\lambda }\times S_{n}$ one may derive a polynomial in $t$
as

$\qquad \qquad f_{\lambda ,\mu }(t)=\sum\limits_{(w,\sigma )\in f_{\lambda
}^{-1}(\mu +\delta (n))}(-t)^{I(\sigma )}$.

\noindent Clearly, one has $f_{\lambda ,\mu }(-1)=Card\{f_{\lambda
}^{-1}(\mu +\delta (n))\}$. On the other hand $f_{\lambda ,\mu
}(1)=K_{\lambda ,\mu }^{-1}$ by Theorem 2.

\bigskip

For two partitions $\mu =(\mu _{1},\ldots ,\mu _{n})$ and $\omega =(\omega
_{1},\ldots ,\omega _{n})$, write $\mu -\omega \in \{r\}$ to simplify the
statement that each of the differences $\mu _{i}-\omega _{i}$ is either $0$
or $1$, and the cardinality of the set $\{i\mid \mu _{i}-\omega _{i}=1\}$ is
precisely $r$ (In the standard terminology of partitions [Ma, p.4-5], $\mu
-\omega \in \{r\}$ is equivalent to the statement that $\omega \subset \mu $
and that the skew diagram $\mu -\omega $ is a vertical $r$-strip). The
notion $\omega <_{r}\mu $ is used to indicate that $\mu ^{(n)}-\omega \in
\{r\}$, where $\mu ^{(n)}=(\mu _{1},\ldots ,\mu _{n-1})$.

For two partitions $\lambda =\{\lambda _{1},\ldots ,\lambda _{k}\},\mu =(\mu
_{1},\ldots ,\mu _{n})\in P(m\mid n)$ denote by $S(\lambda ,\mu )$ the set
of all chains of partitions

$S:\qquad \mu ^{0}=(0)<_{j_{1}}\mu ^{1}<_{j_{2}}\ldots <_{j_{k}}\mu ^{k}=\mu
$

\noindent so that if we assume $\mu ^{i}=(\mu _{1}^{\prime },\ldots ,\mu
_{n}^{\prime })$ and put $b_{i}=\mu _{n}^{\prime }+j_{i}$ for each $1\leq
i\leq k$, then the sequence $(b_{1},\ldots ,b_{k})$ is a permutation of $%
\lambda _{1},\ldots ,\lambda _{k}$.

For a $S\in S(\lambda ,\mu )$ we set $sign(S)=(-1)^{j_{1}+\ldots +j_{k}}$.

\textbf{Theorem 3. }\textsl{If }$\lambda =(r_{1}^{i_{1}},\ldots
,r_{k}^{i_{k}}),\mu =(\mu _{1},\ldots ,\mu _{n})\in P(m\mid n)$\textsl{, then%
}

\noindent (2.1)$\qquad \qquad K_{\lambda ,\mu }^{-1}=\sum\limits_{r_{j}-\mu
_{n}\geq 0}(-1)^{r_{j}-\mu _{n}}\sum\limits_{\mu ^{(n)}-\omega \in
\{r_{j}-\mu _{n}\}}K_{\lambda \lbrack j],\omega }^{-1}$\textsl{.}

\noindent \textsl{Consequently, }$K_{\lambda ,\mu }^{-1}=\sum\limits_{S\in
S(\lambda ,\mu )}sign(S)$\textsl{. }

\bigskip

Usually, a combinatorial identity is achieved whenever a given quantity is
evaluated in two different ways. Combining (1.4) with (2.1) gives

\textbf{Corollary 1}(Identities in the inverse of Kostka matrix). \textsl{%
For any two partitions }$\lambda =(r_{1}^{i_{1}},\ldots ,r_{k}^{i_{k}})$%
\textsl{\ and }$\mu =\{\mu _{1},\ldots ,\mu _{l}\}$\textsl{\ we have}

\noindent $\sum\limits_{\substack{ r_{j}-\mu _{n}\geq 0  \\ \mu
^{(n)}-\omega \in \{r_{j}-\mu _{n}\}}}(-1)^{r_{j}-\mu _{n}}K_{\lambda
\lbrack j],\omega }^{-1}=\sum\limits_{\substack{ 1\leq j\leq k  \\ r_{j}=\mu
_{i}+i-1}}(-1)^{i-1}K_{\lambda \lbrack j],(\mu _{1}-1,\ldots ,\mu
_{i-1}-1,\mu _{i+1},\ldots ,\mu _{l})}^{-1}$\textsl{.}

\bigskip

The non-triviality of these identities can be easily seen from many of their
specializations.

\textbf{Example 1.} Taking $\lambda =(2,3)$ and $\mu =(1^{3},2)$ gives

$\qquad
K_{(3),(1^{3})}^{-1}-K_{(2),(1^{2})}^{-1}=-K_{(3),(1,2)}^{-1}+K_{(2),(2)}^{-1}
$.

\noindent If $\lambda =(1,2^{2})$ and $\mu =(1^{3},2)$ we have

$\qquad
K_{(1,2),(1^{3})}^{-1}=K_{(2^{2}),(1^{2},2)}^{-1}-K_{(1,2),(1,2)}^{-1}$.

\textbf{Remark 2.} In the theory of symmetric functions the \textsl{%
Hall-Littlewood functions} (associated to a partition $\lambda \in P(m\mid
n) $) is defined as

$\qquad P_{\lambda }(x_{1},\ldots ,x_{n};t)=\sum\limits_{w\in S_{n}^{\lambda
}}w(x^{\lambda }\prod\limits_{\lambda _{i}<\lambda _{j}}\frac{x_{j}-tx_{i}}{%
x_{j}-x_{i}})$ [Ma,p.208].

\noindent Since $P_{\lambda }$ is symmetric in $x_{1},\ldots ,x_{n}$ we have
the expression

$\qquad P_{\lambda }(x_{1},\ldots ,x_{n};t)=\sum\limits_{\mu }K_{\lambda
,\mu }^{-1}(t)s_{\mu }$

\noindent for some polynomial $K_{\lambda ,\mu }^{-1}(t)\in \mathbb{Z}[t]$
(cf. [Ma, p.209]) which might be considered as one-parameter variation of $%
K_{\lambda ,\mu }^{-1}$ in view of the fact $K_{\lambda ,\mu
}^{-1}(1)=K_{\lambda ,\mu }^{-1}$, it would be of interest to see if
Corollary 1 can be deduced from certain relations among the polynomials $%
K_{\lambda ,\mu }^{-1}(t)$ (when evaluated at $t=1$).

\bigskip

For $\lambda =(\lambda _{1},\ldots ,\lambda _{n})$, $\mu =(\mu _{1},\ldots
,\mu _{n})\in P(m\mid n)$, write $\lambda >\mu $ to express the fact that
the last non-zero difference $\lambda _{i}-\mu _{i}$ is positive.

\textbf{Corollary 2 (Cancellation principles)}. \textsl{Given }$\lambda
=(\lambda _{1},\ldots ,\lambda _{n})$\textsl{, }$\mu =(\mu _{1},\ldots ,\mu
_{n})\in P(m\mid n)$\textsl{, we have}

\textsl{1) If }$\lambda _{j}=\mu _{j}$\textsl{\ for }$k\leq j\leq n$\textsl{%
, then}

\textsl{\ \qquad }$K_{(\lambda _{1},\ldots ,\lambda _{n}),(\mu _{1},\ldots
,\mu _{n})}^{-1}=K_{(\lambda _{1},\ldots ,\lambda _{k}),(\mu _{1},\ldots
,\mu _{k})}^{-1}$\textsl{.}

\noindent \textsl{In particular, }$K_{\lambda ,\lambda }^{-1}=K_{0,0}^{-1}=1$%
\textsl{\ for all }$\lambda \in P(m\mid n)$\textsl{.}

\textsl{2) If either }$\lambda <\mu $\textsl{\ or }$l(\lambda )>l(\mu ),$
\textsl{then} $K_{\lambda ,\mu }^{-1}=0$\textsl{.}

\textbf{Proof.} If $\mu _{n}\geq \lambda _{n}$, we get immediately from
(2.1) that

\noindent (2.2)$\qquad K_{\lambda ,\mu }^{-1}=\{%
\begin{array}{c}
K_{(\lambda _{1},\ldots ,\lambda _{n-1}),(\mu _{1},\ldots ,\mu _{n-1})}^{-1}%
\text{, if }\mu _{n}=\lambda _{n}\text{;} \\
0\text{, \ if }\mu _{n}>\lambda _{n}\text{.\qquad \qquad \qquad \qquad }%
\end{array}%
$

\noindent This verifies 1) and the first item of 2). The second assertion in
2) is also clear by Theorem 3, since the set $S(\lambda ,\mu )$ must be
empty if $l(\lambda )>l(\mu )$.$\square $

\section{Proofs of Theorem 2 and 3}

Let $\mathbb{Z}[x_{1},\ldots ,x_{n}]$ be the ring of polynomials in $%
x_{1},\ldots ,x_{n}$ and, for an $m\geq 0$, let $\mathbb{Z}[x_{1},\ldots
,x_{n}]^{m}$ be the submodule spanned by the homogeneous polynomials of
degree $m$. By considering $\Lambda _{n}$ a subring of $\mathbb{Z}%
[x_{1},\ldots ,x_{n}]$ in the obvious way, $\Lambda _{n}^{m}$ becomes a
submodule of $\mathbb{Z}[x_{1},\ldots ,x_{n}]^{m}$.

For a sequence $\alpha =(\alpha _{1},\ldots ,\alpha _{n})\in \mathbb{N}^{n}$
consider the additive operator

$\quad \qquad D_{\alpha }:\mathbb{Z}[x_{1},\ldots ,x_{n}]^{\mid \alpha \mid
}\rightarrow \mathbb{Z}$, $\mid \alpha \mid =\Sigma \alpha _{i}$,

\noindent by $D_{\alpha }h=$the coefficient of the monomial $x^{\alpha }$ in
$h$. Alternatively, it can be expressed in terms of partial derivatives as

$\qquad \qquad D_{\alpha }h=\frac{1}{\alpha _{1}!\ldots \alpha _{n}!}\frac{%
\partial ^{\mid \alpha \mid }h}{\partial x_{1}^{\alpha _{1}}\ldots \partial
x_{n}^{\alpha _{n}}}$.

\noindent In [D], [ZD] these operators were applied to express the
integrations along certain flag manifolds and were useful in computing the
degrees of some classical projective varieties.

\bigskip

Since $\{s_{\lambda }(n)\mid \lambda \in P(m\mid n)\}$ is an additive bases
of $\Lambda _{n}^{m}$, every $h\in \Lambda _{n}^{m}$ can be uniquely written
as an integral combination of the $s_{\lambda }(n)$. We give such an
algorithm (cf. [D, Lemma 3.3]).

\textbf{Lemma 1.} If $h=\underset{\lambda \in P(m\mid n)}{\Sigma }c_{\lambda
}s_{\lambda }(n)$, then $c_{\lambda }=D_{\lambda +\delta (n)}(ha_{\delta
(n)})$.

\textbf{Proof}. For a $\mu =(r_{1},\ldots ,r_{n})\in P(m\mid n)$ multiplying
both sides of $h=\underset{\lambda \in P(m\mid n)}{\Sigma }c_{\lambda
}s_{\lambda }$ by $a_{\delta (n)}$, and applying $D_{\mu +\delta (n)}$ yields

$\qquad \qquad D_{\mu +\delta (n)}(ha_{\delta (n)})=\underset{\lambda \in
P(m\mid n)}{\Sigma }c_{\lambda }D_{\mu +\delta (n)}(a_{\lambda +\delta (n)})$%
.

\noindent The proof is completed by observing that the monomial $x^{\mu
+\delta (n)}$ is a term in $a_{\lambda +\delta (n)}$ if and only if $\lambda
=\mu $ (note that the sequences $\lambda +\delta (n)$ and $\mu +\delta (n)$
are strict increasing), and its coefficient in $a_{\mu +\delta (n)}$ is
precisely $1.\square $

\bigskip

Theorem 1 follows directly from Lemma 1.

\textbf{Proof of Theorem 1. }From Lemma 1 we have $\ \ \ \ $

$\qquad K_{\lambda ,\mu }=D_{\mu +\delta (n)}(m_{\lambda }(n)a_{\delta (n)})$

$=D_{\mu +\delta (n)}(\sum\limits_{w\in S_{n}^{\lambda }}x^{w(\lambda
)}\sum\limits_{\sigma \in S_{n}}\varepsilon (\sigma )x^{\sigma (\delta
(n))})=\sum\limits_{\substack{ (w,\sigma )\in S_{n}^{\lambda }\times S_{n}
\\ w(\lambda )+\sigma (\delta (n))=\mu +\delta (n)}}\varepsilon (\sigma )$.$%
\square $

\bigskip

The proof of Theorem 2 requires a little more preparation. The first of
these is the following ``elimination law''(cf. [D, Section 3]).

\textbf{Lemma 2}. \textsl{If }$h(x_{1},\ldots ,x_{n})=\Sigma h_{i}x_{n}^{i}$%
\textsl{, with }$h_{i}\in Z[x_{1},\ldots ,x_{n-1}]$\textsl{, then}

$\QTR{sl}{\qquad D}_{(r_{1},\ldots ,r_{n})}\QTR{sl}{h(x}_{1}\QTR{sl}{,}$%
\textsl{$\ldots $}$\QTR{sl}{,x}_{n}\QTR{sl}{)=D}_{(r_{1},\ldots ,r_{n-1})}%
\QTR{sl}{h}_{r_{n}}\QTR{sl}{(x}_{1}\QTR{sl}{,}$\textsl{$\ldots $}$\QTR{sl}{,x%
}_{n-1}\QTR{sl}{).}$

\bigskip

Let $e_{r}(n)\in \Lambda _{n}^{r}$\ be the $r^{th}$\ elementary symmetric
function, $0\leq r\leq n$.\textsl{\ }Then $e_{r}(n)=m_{(1^{r})}(n)$. The
next result (cf. [Ma, p.73]) is known as the \textsl{Pieri-formula}.

\textbf{Lemma 3.} \textsl{\ For any }$\lambda \in P(m\mid n)$\textsl{, }$%
s_{\lambda }(n)e_{r}(n)=\sum\limits_{\mu -\lambda \in \{r\}}s_{\mu }(n)$%
\textsl{.}

\bigskip

The correspondence $P(m\mid n-1)\rightarrow P(m\mid n)$ by $\lambda
=(\lambda _{1},\ldots ,\lambda _{n-1})\rightarrow \lambda ^{\prime
}=(0,\lambda _{1},\ldots ,\lambda _{n-1})$ is a bijection whenever $m\leq
n-1 $ [Ma, p.18]. The stability of the numbers $K_{\lambda ,\mu }^{-1}$ can
now be stated as

\textbf{Lemma 4.} \textsl{If }$\lambda ,\mu \in P(m\mid n-1)$\textsl{\ and
if }$m\leq n-1$\textsl{, then }$K_{\lambda ,\mu }^{-1}=K_{\lambda ^{\prime
},\mu ^{\prime }}^{-1}$\textsl{.}

\noindent In view of this we do not differentiate between $\lambda $ and $%
\lambda ^{\prime }$.

\bigskip

For a partition $\lambda =(\lambda _{1},\ldots ,\lambda _{n})\in P(m\mid n)$
and for $1\leq i\leq n$, write $\lambda ^{(i)}$ for the partition $(\lambda
_{1},\ldots ,\lambda _{i-1},\lambda _{i+1}+1,\ldots ,\lambda _{n}+1)\in
P(m^{\prime }\mid n-1)$, where $m^{\prime }=m-\lambda _{i}+(n-i)$.

\textbf{Proof of Theorem 2. }Suppose that for all $m\leq n$ and $\lambda
,\mu \in P(m\mid n)$ the numbers $K_{\lambda ,\mu }^{-1}(n)\in \mathbb{Z}$
are defined by

\noindent (3.1) $\qquad m_{\lambda }(n)=\sum\limits_{\mu \in P(m\mid
n)}K_{\lambda ,\mu }^{-1}(n)s_{\mu }(n)$

With $\lambda =(r_{1}^{i_{1}},\ldots ,r_{k}^{i_{k}})$, the formula
expressing $m_{\lambda }(n)$ as a polynomial in $x_{n}$ with coefficients
from $\mathbb{Z}[x_{1},\ldots ,x_{n-1}]$ is

$\qquad m_{\lambda }(n)=m_{\lambda }(n-1)+\sum\limits_{1\leq j\leq
k}m_{\lambda \lbrack j]}(n-1)x_{n}^{r_{j}}$.

\noindent Expanding the determinant $a_{\delta (n)}$ with respect to the
last column yields

\qquad $a_{\delta (n)}=\sum\limits_{1\leq i\leq n}(-1)^{i-1}a_{\delta
(n)^{(n-i+1)}}x_{n}^{n-i}$.

\noindent It follows that

$m_{\lambda }(n)a_{\delta (n)}=m_{\lambda }(n-1)a_{\delta (n)}+\sum\limits
_{\substack{ 1\leq i\leq n  \\ 1\leq j\leq k}}(-1)^{i-1}m_{\lambda \lbrack
j]}(n-1)a_{\delta (n)^{(n-i+1)}}x_{n}^{r_{j}+n-i}$.

>From Lemma 1 we have

\noindent (3.2) $\qquad K_{\lambda ,\mu }^{-1}(n)=D_{\mu +\delta
(n)}(m_{\lambda }(n)a_{\delta (n)})$.

\noindent Since we can assume $\mu _{n}\geq 1$, the coefficient of $%
x_{n}^{\mu _{n}+n-1}$ in $m_{\lambda }(n-1)a_{\delta (n)}$ is clearly zero.
Applying the elimination law (Lemma 2) to the right hand side of (3.2) gives

\noindent (3.3)$K_{\lambda ,\mu }^{-1}(n)=D_{\mu ^{(n)}+\delta
(n-1)}(\sum\limits_{\substack{ \mu _{n}+n-1=r_{j}+n-i  \\ 1\leq j\leq k}}%
(-1)^{i-1}m_{\lambda \lbrack j]}(n-1)a_{\delta (n)^{(n-i+1)}})$

$\quad =D_{\mu ^{(n)}+\delta (n-1)}\{[\sum\limits_{\substack{ i=r_{j}+1-\mu
_{n}  \\ 1\leq j\leq k}}(-1)^{i-1}m_{\lambda \lbrack
j]}(n-1)e_{i-1}(n-1)]a_{\delta (n-1)}\}$,

\noindent where $\mu ^{(n)}=(\mu _{1},\ldots ,\mu _{n-1})$($\in P(m-\mu
_{n}\mid n-1)$), and where the second equality follows from $a_{\delta
(n)^{(n-i+1)}}=e_{i-1}(n-1)a_{\delta (n-1)}$. Since

$\qquad m_{\lambda \lbrack j]}(n-1)=\sum\limits_{\omega \in P(m-r_{j}\mid
n-1)}K_{\lambda \lbrack j],\omega }^{-1}(n-1)s_{\omega }(n-1)$

\noindent by our assumption (3.1) and since

\ $\qquad s_{\omega }(n-1)e_{i-1}(n-1)=\sum\limits_{\gamma -\omega \in
\{i-1\}}s_{\gamma }(n-1)$

\noindent by Lemma 3, we get from (3.3) that

\ $\qquad K_{\lambda ,\mu }^{-1}(n)=\sum\limits_{r_{j}-\mu _{n}\geq
0}(-1)^{r_{j}-\mu _{n}}\sum\limits_{\mu ^{(n)}-\omega \in \{r_{j}-\mu
_{n}\}}K_{\lambda \lbrack j],\omega }^{-1}(n-1)$

\noindent (again by Lemma 1). Finally, $\lambda \lbrack j],\omega \in $ $%
P(m-r_{j}\mid n-1)$ implies that

$\qquad K_{\lambda ,\mu }^{-1}(n)=\sum\limits_{r_{j}-\mu _{n}\geq
0}(-1)^{r_{j}-\mu _{n}}\sum\limits_{\mu ^{(n)}-\omega \in \{r_{j}-\mu
_{n}\}}K_{\lambda \lbrack j],\omega }^{-1}(n)$

\noindent by Lemma 4. This completes the proof.$\square $

\section{Applications}

Theorem 3 (i.e. (2.1)) enables one to deduce explicit formulas for $%
K_{\lambda ,\mu }^{-1}$ for special cases of $\lambda $ and $\mu $. This
section is devoted to such examples that have interesting numerical features.

\bigskip

Consider firstly $\lambda =(r_{1}^{i_{1}},\ldots ,r_{k}^{i_{k}})$, $\mu
=(1^{m})\in P(m\mid n)$. We have

\noindent (4.1)$\qquad \qquad K_{\lambda ,(1^{m})}^{-1}=\sum\limits_{1\leq
j\leq k}(-1)^{r_{j}-1}K_{\lambda \lbrack j],(1^{m-r_{j}})}^{-1}$.

\noindent by (2.1). An induction on $l(\lambda )=\Sigma i_{j}$ verifies (cf.
[ER, Corollary 1]):

\textbf{Lemma 5.} $K_{\lambda ,(1^{m})}^{-1}=\frac{(-1)^{l(\mu )-l(\lambda
)}l(\lambda )!}{i_{1}!i_{2}!\ldots i_{k}!}$.

\textbf{Proof.} $K_{(m),(1^{m})}^{-1}=(-1)^{m-1}$ by Theorem 3. From the
inductive hypothesis one has

$\qquad K_{\lambda \lbrack j],(1^{m-r_{j}})}^{-1}=\frac{(-1)^{m-r_{j}-l(%
\lambda )-1}(l(\lambda )-1)!}{i_{1}!\ldots i_{j-1}!(i_{j}-1)!i_{j+1}!\ldots
i_{k}!}$, $1\leq j\leq k$.

\noindent Substituting these in (4.1) yields

$\qquad K_{\lambda ,(1^{m})}^{-1}=\frac{(-1)^{m-l(\lambda )}(l(\lambda )-1)!%
}{i_{1}!i_{2}!\ldots i_{k}!}\Sigma i_{j}=\frac{(-1)^{l(\mu )-l(\lambda
)}l(\lambda )!}{i_{1}!i_{2}!\ldots i_{k}!}$.$\square $

\bigskip

For a $\lambda \in P(m\mid n)$ we set $\lambda _{i}^{\prime }=Card\{j\mid
\lambda _{j}\geq i\}$. Continuing from Lemma 5 we recover the following
result originally due to Kostka [K2]

\textbf{Lemma 6.} $K_{\lambda ,(1^{m}a)}^{-1}=\frac{(-1)^{l(\mu )-l(\lambda
)}(l(\lambda )-1)!}{i_{1}!i_{2}!\ldots i_{k}!}\lambda _{a}^{\prime }$.

\textbf{Proof.} $K_{\lambda ,(1^{m}a)}^{-1}=$ $\sum\limits_{r_{j}\geq
a}(-1)^{r_{j}-a}K_{\lambda \lbrack j],(1^{m-r_{j}+a})}^{-1}$ (by (2.1))

$\qquad =\sum\limits_{r_{j}\geq a}(-1)^{r_{j}-a}\frac{(-1)^{m-r_{j}+a-l(%
\lambda )+1}(l(\lambda )-1)!}{i_{1}!\ldots i_{j-1}!(i_{j}-1)!i_{j+1}!\ldots
i_{k}!}$ (by Lemma 5)

$\qquad =\frac{(-1)^{l(\mu )-l(\lambda )}(l(\lambda )-1)!}{%
i_{1}!i_{2}!\ldots i_{k}!}\lambda _{a}^{\prime }$.$\square $

\bigskip

An extension of Lemma 6 is the case $\mu =(1^{m},a,b)$, $1<a\leq b$. We may
assume that $\lambda =(r_{1}^{i_{1}},\ldots ,r_{k}^{i_{k}})\geq \mu $ by 2)
of Corollary 2. Then there exists a unique $d\leq k$ such that $r_{d}\geq b$
but $r_{i}<b$ for all $1\leq i<d$. Consider (for instance) the case $r_{d}>b$%
. We have

\noindent $K_{\lambda ,(1^{m},a,b)}^{-1}=\sum\limits_{d\leq j\leq
k}(-1)^{r_{j}-b}[K_{\lambda \lbrack j],(1^{m-r_{j}+b},a)}^{-1}+K_{\lambda
\lbrack j],(1^{m-r_{j}+b+1},a-1)}^{-1}]$ (by (2.1))

\noindent $=\sum\limits_{d\leq j\leq k}(-1)^{r_{j}-b}[\frac{%
(-1)^{m-r_{j}+b-l(\lambda )+1}(l(\lambda )-2)!}{i_{1}!\ldots
i_{j-1}!(i_{j}-1)!i_{j+1}!\ldots i_{k}!}(\lambda \lbrack j]_{a}^{\prime
}-\lambda \lbrack j]_{a-1}^{\prime })]$( by Lemma 6 )

\noindent $=\frac{(-1)^{l(\mu )-l(\lambda )}(l(\lambda )-2)!}{%
i_{1}!i_{2}!\ldots i_{k}!}\sum\limits_{d\leq j\leq k}i_{j}(\lambda \lbrack
j]_{a-1}^{\prime }-\lambda \lbrack j]_{a}^{\prime })$.

\noindent Thus, for a $\omega =(\omega _{1},\ldots ,\omega _{k})$, if we let
$D_{c}(\omega )=Card\{i\mid \omega _{i}=c\}$, we have

\textbf{Corollary 3.} \textsl{If }$r_{d}>b$\textsl{, then}

$\qquad K_{\lambda ,(1^{m},a,b)}^{-1}=\frac{(-1)^{l(\mu )-l(\lambda
)}(l(\lambda )-2)!}{i_{1}!i_{2}!\ldots i_{k}!}\sum\limits_{d\leq j\leq
k}i_{j}D_{a-1}(\lambda \lbrack j])$\textsl{;}

\noindent \textsl{If }$r_{d}=b$\textsl{, then}

$\qquad K_{\lambda ,(1^{m},a,b)}^{-1}=\frac{(-1)^{l(\mu )-l(\lambda
)}(l(\lambda )-2)!}{i_{1}!i_{2}!\ldots i_{k}!}[i_{d}\lambda \lbrack
d]_{a}^{\prime }+\sum\limits_{d+1\leq j\leq k}i_{j}D_{a-1}(\lambda \lbrack
j])]$\textsl{.}

\bigskip

Let us compute $K_{(1^{k},2^{l}),\mu }^{-1}$. Firstly $K_{(1^{k},2^{l}),\mu
}^{-1}=0$ if $\mu \neq (1^{k+2t},2^{l-t})$ by 2) of Corollary 2. Moreover

$\qquad
K_{(1^{k},2^{l}),(1^{k+2t},2^{l-t})}^{-1}=K_{(1^{k},2^{t}),(1^{k+2t})}^{-1}$
(by 1) of Corollary 2)

$\qquad \qquad =\frac{(-1)^{t}(k+t)!}{k!t!}$ (By Lemma 5).

\noindent Summarizing we get (cf. also [ER, Corollary 3])

\textbf{Corollary 4.} $m_{(1^{k},2^{l})}(n)=\sum\limits_{0\leq t\leq l}\frac{%
(-1)^{t}(k+t)!}{k!t!}s_{(1^{k+2t},2^{l-t})}(n)$\textsl{.}

\textbf{Remark 3.} Combining Corollary 4 with the Giambelli-formula

$\qquad s_{(1^{m-k},2^{k})}(n)=e_{k}(n)e_{m}(n)-e_{k-1}(n)e_{m+1}(n)$

\noindent yields an integral lift of the Wu-formula (1.3).$\square $

\bigskip

Let us compute $K_{(1^{k},3^{l}),\mu }^{-1}$. Again, By 2) of Corollary 2 we
have $K_{(1^{k},3^{l}),\mu }^{-1}=0$ unless $\mu =(1^{a},2^{b},3^{c})$ with $%
c\leq l$. Furthermore, 1) of Corollary 2 implies that $%
K_{(1^{k},3^{l}),(1^{a},2^{b},3^{c})}^{-1}=K_{(1^{k},3^{l-c}),(1^{a},2^{b})}^{-1}
$. So it remains to find $K_{(1^{k},3^{l}),(1^{a},2^{b})}^{-1}$, where $%
k+3l=a+2b$.

For fixed $k$ and $l$ consider the polynomial in $t$

\noindent (4.2)$\qquad \qquad g_{k,l}(t)=\sum\limits_{0\leq
b,k+3l=a+2b}K_{(1^{k},3^{l}),(1^{a},2^{b})}^{-1}t^{b}$.

\noindent Since

$K_{(1^{k},3^{l}),(1^{a},2^{b})}^{-1}=\{%
\begin{array}{c}
\frac{(k+l)!}{k!l!}\text{ if }b=0\text{, }a\geq 1\text{ (Lemma 5);\quad
\qquad \qquad } \\
-K_{(1^{k},3^{l-1}),(1^{a-1})}^{-1}\text{ if }b=1\text{,}a\geq 1\text{ (by
(2.1));} \\
-K_{(1^{k},3^{l-1}),(1,2^{b-2})}^{-1}\text{ if }a=0\text{ (by (2.1)),\quad
\quad }%
\end{array}%
$

\noindent and since

$%
K_{(1^{k},3^{l}),(1^{a},2^{b})}^{-1}=-K_{(1^{k},3^{l-1}),(1^{a-1},2^{b-1})}^{-1}-K_{(1^{k},3^{l-1}),(1^{a+1},2^{b-2})}^{-1}
$

\noindent if $b\geq 2$, $a\geq 1$ (by (2.1)), we obtain

\noindent (4.3)$\qquad \qquad g_{k,l}(t)=\frac{(k+l)!}{k!l!}%
-(t+t^{2})g_{k,l-1}(t)+\varepsilon (t)$,

\noindent where

$\qquad \varepsilon (t)=\{%
\begin{array}{c}
0\text{ \qquad if }k+3l\text{ is even;\qquad \qquad \qquad \qquad \qquad }
\\
K_{(1^{k},3^{l}),(2^{\frac{k+3(l-1)}{2}})}^{-1}t^{\frac{k+3(l-1)}{2}}\text{
\qquad if }k+3l\text{ is odd.}%
\end{array}%
$

\bigskip

Assume now that $k>l-1$. Then $\varepsilon (t)=0$ by 2) of Corollary 2. We
infer from (4.3) and the obvious relation $g_{k,0}(t)=(-1)^{k-1}$ (by Lemma
5) that

\textbf{Corollary 5. }\textsl{If }$k>l-1$\textsl{, then }$%
g_{k,l}(t)=\sum\limits_{0\leq i\leq l}(-1)^{i}\frac{(k+l-i)!}{k!(l-i)!}%
(t+t^{2})^{i}$\textsl{.}

\bigskip

In order to apply (4.3) to find an expression of $g_{k,l}(t)$ for the case $%
k\leq l-1$, we need to compute $K_{(1^{k},3^{l}),(2^{b})}^{-1}$ (for $k+3l$
even). This will be done by combining (2.1) with\textbf{\ }the
Egecioglu-Remmel formula (1.4). Consider, for a fixed $b$, the polynomial in
$t$

$\qquad \qquad
h_{b}(t)=\sum\limits_{k+3l=2b}K_{(1^{k},3^{l}),(2^{b})}^{-1}t^{k}$.

\noindent From

$K_{(1^{k},3^{l}),(2^{b})}^{-1}=-K_{(1^{k},3^{l-1}),(1,2^{b-2})}^{-1}$ (by
(2.1))

$\qquad =\{%
\begin{array}{c}
-K_{(1^{k-1},3^{l-1}),(2^{b-2})}^{-1}+K_{(1^{k},3^{l-2}),(2^{b-3})}^{-1}%
\text{ if }k\geq 1\text{;} \\
K_{(3^{l-2}),(2^{b-3})}^{-1}\text{ if }k=0\qquad \text{\qquad \qquad \qquad
\qquad }%
\end{array}%
$

\noindent (by (1.4)) we get

$\qquad h_{b}(t)=-th_{b-2}(t)+h_{b-3}(t)$.

\noindent It follows that

$\qquad h_{b}(t)=(t^{2},-2t,1)\left(
\begin{array}{c}
h_{b-4} \\
h_{b-5} \\
h_{b-6}%
\end{array}%
\right) $

\noindent and that

$\qquad \left(
\begin{array}{c}
h_{b-4} \\
h_{b-5} \\
h_{b-6}%
\end{array}%
\right) =A(t)\left(
\begin{array}{c}
h_{b-6} \\
h_{b-7} \\
h_{b-8}%
\end{array}%
\right) $, where $A(t)=\left(
\begin{array}{ccc}
-t & 1 & 0 \\
0 & -t & 1 \\
1 & 0 & 0%
\end{array}%
\right) $.

\textbf{Corollary 6.} \textsl{Assume that }$b=2k+r$\textsl{\ with }$0\leq
r\leq 1$\textsl{. Then}

$\qquad h_{b}(t)=(t^{2},-2t,1)A(t)^{k-3}\left(
\begin{array}{c}
h_{2+r}(t) \\
h_{1+r}(t) \\
h_{r}(t)%
\end{array}%
\right) $\textsl{, }

\noindent \textsl{where }$(h_{0}(t),h_{1}(t),h_{2}(t),h_{3}(t))=(1,0,-t,1)$%
\textsl{.}

The last equation can be easily obtained by using either (1.4) or (2.1).

\bigskip

\textbf{Example 2. }Based on Corollary 6 a program to expand $h_{b}(t)$ has
been compiled. We list below the results for $25\leq b\leq 30$ produced by
the program.

$\qquad h_{25}(t)=36t^{2}-252t^{5}+165t^{8}-12t^{11};$

$\qquad h_{26}(t)=-9t+210t^{4}-330t^{7}+66t^{10}-t^{13};$

$\qquad h_{27}(t)=1-120t^{3}+462t^{6}-220t^{9}+13t^{12};$

$\qquad h_{28}(t)=45t^{2}-462t^{5}+495t^{8}-78t^{11}+t^{14};$

$\qquad h_{29}(t)=-10t+330t^{4}-792t^{7}+286t^{10}-14t^{13};$

$\qquad h_{30}(t)=1-165t^{3}+924t^{6}-715t^{9}+91t^{12}-t^{15}$.

By the discussion in Section 1 on the Steenrod cohomology operations, the
above polynomials can be used to reveal deep information on the topology of
the classifying space $BU(n)$ (as well as the complex Grassmannians). For
instance, we see from

$\qquad h_{30}(t)\equiv 1+2t^{9}+t^{12}+2t^{15}$ (mod $3$)

\noindent that the attaching map of the Schubert cell $\Omega (2^{30})$,
which has the real dimension $120$, cannot avoid any of the Schubert cells $%
\Omega (1^{10})$, $\Omega (1^{16})$, $\Omega (1^{18})$ and $\Omega (1^{20})$
via homotopies.

\bigskip

\textbf{Acknowledgements.} The author feels very grateful to his referee for
many improvements on an earlier version of this paper. Thanks are also due
to Dr. Xuezhi Zhao who wrote a program to implement the polynomial $h_{b}(t)$
in Corollary 6 (cf. Example 2).

\begin{center}
\textbf{References}
\end{center}

[D] H. Duan, Some enumerative formulas for flag manifolds, Communications in
Algebra, 29(10), (2001), 4395-4419.

[K1] C. Kostka, Uber den Zusammenhang zwischen einigen Formen von
symmetrischen Funktionen, J. Reine Angew. Math., 93(1882), 89-123.

[K2] C. Kostka, Tafeln fur symmetrische Funktionen bis zur elften Dimension,
Wissenschaftliche Beilage zum Pragramm des konigl. Gymnasiums und
Realgymnasiums zu Insterberg, 1908.

[ER] O. Egecioglu and J. B. Remmel, A combinatorial interpretation of the
inverse Kostka matrix, Linear and Multilinear algebra, 26(1990), 59-84.

[L] C. Lenart, The combinatorial of Steenrod operations on the cohomology of
Grassmannians, Advances in Math. 136(1998), 251-283.

[Ma] I. G. Macdonald, Symmetric functions and Hall polynomials, Oxford
Mathematical Monographs, Oxford University Press, Oxford, second ed., 1995.

[MS] J. Milnor and J. Stasheff, Characteristic classes, Ann. of Math.
Studies 76, Princeton Univ. Press, 1975.

[S] P. Shay, mod-p Wu formulas for the Steenrod algebra and the Dyer-Lashof
algebra, Proc. AMS, 63(1977), 339-347.

[St] N. E. Steenrod and D. B. A. Epstein, Cohomology Operations, Ann. of
Math. Stud., Princeton Univ. Press, Princeton, NJ, 1962.

[W] W.T.Wu, Les i-carries dans une variete grassmanniene, C. R. Acad. Sci.
Paris 230(1950), 918-920).

[ZD] X. Zhao and H. Duan, A Mathematica Program for the Degrees of certain
Schubert varieties, J. Symbolic Computation, 33(2002), 507-517.

\end{document}